\documentclass[12pt]{article}
\usepackage{amssymb,graphicx}
\usepackage[only,fivrm]{rawfonts}
\usepackage{t1enc}
\input{prepictex}
\input{pictex}
\input{postpictex}

\newcommand{\wt}{\widetilde}
\newcommand{\R}{\mathbb R}
\newcommand{\C}{\mathbb C}
\newcommand{\Q}{\mathbb Q}

\begin{document}
\thispagestyle{empty}
\footnotetext{
\par
\noindent
{\bf Abbreviated title:} Beckman-Quarles type theorems.}
\hfill
\vskip 1.2truecm
\centerline{{\large Beckman-Quarles type theorems for mappings from ${\R}^n$ to ${\C}^n$}}
\vskip 1.2truecm
\centerline{{\large \sc Apoloniusz Tyszka}}
\vskip 1.2truecm
\par
\noindent
{\bf Summary.} Let $\varphi_n:{\C}^n \times {\C}^n \to \C$,
$\varphi_n((x_1,...,x_n),(y_1,...,y_n))=(x_1-y_1)^2+...+(x_n-y_n)^2$.
We say that $f:{\R}^n \to {\C}^n$ preserves distance $d \geq 0$ if
for each $x,y \in {\R}^n$ $\varphi_n(x,y)=d^2$ implies
$\varphi_n(f(x),f(y))=d^2$. Let $A_n$ denote the set of all
positive numbers~$d$ such that any map $f:{\R}^n \to {\C}^n$
that preserves unit distance preserves also distance~$d$.
Let $D_{n}$ denote the set of all positive numbers $d$
with the property:
if~$x,y \in {\R}^n$ and $|x-y|=d$ then there exists a finite set $S_{xy}$
with $\left
\{x,y \right\} \subseteq S_{xy} \subseteq {\R}^n$ such that any map
$f:S_{xy}\rightarrow {\C}^n$ that preserves unit distance
preserves also the distance between $x$ and $y$.
Obviously, $D_n \subseteq A_n$.
We prove:
{\bf (1)} $A_n \subseteq \left\{d>0: d^2 \in \Q \right\}$,
{\bf (2)} for $n \geq 2$ $D_n$ is a dense subset of $(0,\infty)$.
Item {\bf (2)} implies that each mapping $f$
from ${\R}^n$ to ${\C}^n$ ($n \geq 2$)
preserving unit distance preserves all distances,
if $f$ is continuous with respect to the product topologies
on ${\R}^n$ and ${\C}^n$.
\vskip 1.2truecm
\par
The classical Beckman-Quarles theorem states that each unit-distance
preserving mapping from ${\R}^n$ to ${\R}^n$ ($n \geq 2$) is an isometry,
see \cite{Beckman}, \cite{Benz}, \cite{Everling} and~\cite{Lester}.
By a {\sl complex isometry} of ${\C}^n$ we understand any map
$f:{\C}^n \rightarrow {\C}^n$ of the form
$$
f(z_1,z_2,...,z_n)=(z_1',z_2',...,z_n')
$$
where
$$
z_j'=a_{0j}+a_{1j}z_{1}+a_{2j}z_{2}+...+a_{nj}z_{n} \hspace{0.5truecm}(j=1,2,...,n),
$$
the coefficients $a_{ij}$ are complex and the matrix $||a_{ij}||$
$(i,j=1,2,...,n)$ is orthogonal i.e. satisfies the condition
$$
\sum_{j=1}^{n} a_{\mu j} a_{\nu j}=\delta_{\nu}^{\mu} \hspace{0.5truecm} (\mu,\nu=1,2,...,n)
$$
with \@ Kronecker's \@ delta. \@ Let \@
$\varphi_n:{\C}^n \times {\C}^n \to \C$, \@
$\varphi_n((x_1,...,x_n),(y_1,...,y_n))=(x_1-y_1)^2+...+(x_n-y_n)^2$.
The function $\varphi_n$ is not the square of a distance
because it admits values that are not real.
Nevertheless, according to \cite{Borsuk}
$\varphi_n(x,y)$ is invariant under complex isometries i.e. for every
complex isometry $f:{\C}^n \rightarrow {\C}^n$
$$
\forall x,y \in {\C}^n ~~\varphi_n(f(x),f(y))=\varphi_n(x,y).
$$
Therefore,
we say that $f:{\R}^n \to {\C}^n$ preserves distance $d \geq 0$
if for each $x,y \in {\R}^n$ $\varphi_n(x,y)=d^2$ implies
$\varphi_n(f(x),f(y))=d^2$.
\vskip 0.3truecm
\par
\noindent
By a field endomorphism of $\C$ we understand any map
$g:\C \to \C$ satisfying:
\par
\centerline
{$\forall x,y \in {\C}~~g(x+y)=g(x)+g(y)$,}
\centerline{$\forall x,y \in {\C}~~g(x \cdot y)=g(x) \cdot g(y)$,}
\centerline{$g(0)=0$,}
\centerline{$g(1)=1$.}
\par
\noindent
Bijective endomorphisms are called automorphisms,
for more information on field endomorphisms and automorphisms
of $\C$ the reader is referred to \cite{Kuczma} and \cite{Yale}.
If $r$ is a rational number, then $g(r)=r$ for any field
endomorphism $g:\C \to \C$.
Proposition 1 shows that only rational numbers $r$ have this property:
\vskip 0.3truecm
\par
\noindent
{\bf Proposition 1.} If $r \in \C$ and $r$ is not a rational number,
then there exists a field automorphism $g:\C \to \C$ such that
$g(r) \neq r$.
\vskip 0.3truecm
\par
\noindent
{\it Proof.} Note first that if $E$ is any subfield of $\C$
and if $g$ is an automorphism of $E$, then~$g$ can be extended
to an automorphism of $\C$. This follows from [5, corollaire 1 on p.~109]
by taking $\Omega=\C$ and $K=\Q$.
Now let $r \in \C \setminus \Q$.
If $r$ is algebraic over $\Q$, let $E$ be the
splitting field in $\C$ of the minimal polynomial $\mu$ of $r$ over $\Q$
and let $r' \in E$ be any other root of $\mu$.
Then there exists an automorphism $g$ of $E$ that sends $r$ to $r'$,
see for example [7, corollary 2 on p.~66].
If $r$ is transcendental over $\Q$, let $E=\Q(r)$ and let
$r' \in E$ be any other generator of $E$ (e.g. $r'=1/r$). Then there
exists an automorphism~$g$ of $E$ that sends $r$ to $r'$.
In each case, $g$ can be extended to an automorphism of $\C$.
\vskip 0.3truecm
\par
If $g:\C \to \C$ is a field endomorphism then
$(g_{|\R},...,g_{|\R}): {\R}^n \to {\C}^n$ preserves all distances
$\sqrt{r}$ with rational $r \geq 0$.
Indeed, if $(x_1-y_1)^2+...+(x_n-y_n)^2=(\sqrt{r})^2$ then
\newpage
\begin{eqnarray*}
\varphi_n((g_{|\R},...,g_{|\R})(x_1,...,x_n),(g_{|\R},...,g_{|\R})(y_1,...,y_n))&=&\\
\varphi_n((g(x_1),...,g(x_n)),(g(y_1),...,g(y_n)))&=&\\
(g(x_1)-g(y_1))^2+...+(g(x_n)-g(y_n))^2&=&\\
(g(x_1-y_1))^2+...+(g(x_n-y_n))^2&=&\\
g((x_1-y_1)^2)+...+g((x_n-y_n)^2)&=&\\
g((x_1-y_1)^2+...+(x_n-y_n)^2)&=&g((\sqrt{r})^2)=g(r)=r=(\sqrt{r})^2.
\end{eqnarray*}
\vskip 0.3truecm
\par
\noindent
Theorem 1 shows that we cannot expect a result like
that of Beckman and Quarles:
\vskip 0.3truecm
\par
\noindent
{\bf Theorem~1.} If $x,y \in {\R}^n$ and $|x-y|^2$ is not a rational
number, then there exists $f:{\R}^n \to {\C}^n$ that does not preserve the
distance between $x$ and~$y$ although $f$ preserves all distances
$\sqrt{r}$ with rational $r \geq 0$.
\vskip 0.3truecm
\par
\noindent
{\it Proof.} There exists an isometry $I:{\R}^n \to {\R}^n$
such that $I(x)=(0,0,...,0)$ and $I(y)=(|x-y|,0,...,0)$.
By Proposition 1 there exists
a field automorphism $g:\C \to \C$ such that
$g(|x-y|^2) \neq |x-y|^2$. Thus $g(|x-y|) \neq |x-y|$ and
$g(|x-y|) \neq -|x-y|$. Therefore
$(g_{|\R},...,g_{|\R}): {\R}^n \to {\C}^n$
does not preserve the distance between
$(0,0,...,0) \in {\R}^n$ and $(|x-y|,0,...,0) \in {\R}^n$
although $(g_{|\R},...,g_{|\R})$ preserves all
distances $\sqrt{r}$ with rational $r \geq 0$.
Hence $f:=(g_{|\R},...,g_{|\R}) \circ I: {\R}^n \to {\C}^n$
does not preserve the distance between $x$ and $y$ although
$f$ preserves all distances $\sqrt{r}$ with rational $ r \geq 0$.
\vskip 0.3truecm
\par
Let $A_{n}$ denote the set of all positive numbers $d$
such that any map $f:{\R}^n \to {\C}^n$ that preserves
unit distance preserves also distance $d$.
By Theorem 1 $A_n \subseteq \left\{d>0: d^2 \in \Q \right\}$.
Let $D_{n}$ denote the set of all positive numbers $d$
with the following property:
\vskip 0.3truecm
\par
\noindent
if $x,y \in {\R}^n$ and $|x-y|=d$ then there exists a finite set $S_{xy}$
with $\left\{x,y \right\} \subseteq S_{xy} \subseteq {\R}^n$ such that any map
$f:S_{xy}\rightarrow {\C}^n$ that preserves unit distance
preserves also the distance between $x$ and $y$.
\vskip 0.3truecm
\par
\noindent
In case of $|x-y|=1$ we define $S_{xy}:=\left\{x,y\right\}$. Thus, $1 \in D_n$.
Obviously, $D_n \subseteq A_n$.
\vskip 0.3truecm
\par
We shall study the set $D_n$ and a unit-distance
preserving mapping $f$ from ${\R}^n$ to~${\C}^n$.
We need the following technical Propositions 2-7.
\newpage
\par
\noindent
{\bf Proposition~2} (cf. \cite{Blumenthal}, \cite{Borsuk}).
The points $c_{1}=(z_{1,1},...,z_{1,n}),...,
c_{n+1}=(z_{n+1,1},...,z_{n+1,n}) \in {\C}^n$ are affinely dependent
if and only if their Cayley-Menger determinant
$\Delta(c_1,...,c_{n+1}):=$
$$
\det \left[
\begin{array}{ccccc}
 0  &  1                       &  1                       & ... & 1                       \\
 1  &  0                       & \varphi_n(c_{1},c_{2})   & ... & \varphi_n(c_{1},c_{n+1})\\
 1  & \varphi_n(c_{2},c_{1})   &  0                       & ... & \varphi_n(c_{2},c_{n+1})\\
... & ...                      & ...  	                  & ... & ...                     \\
 1  & \varphi_n(c_{n+1},c_{1}) & \varphi_n(c_{n+1},c_{2}) & ... & 0                       \\
\end{array}\;\right]
$$
\par
\noindent
equals $0$.
\vskip 0.3truecm
\par
\noindent
{\it Proof.} It follows from the equality
$$
\left(
\det \left[
\begin{array}{ccccc}
z_{1,1}   & z_{1,2}   & ... &  z_{1,n}  & 1  \\
z_{2,1}   & z_{2,2}   & ... &  z_{2,n}  & 1  \\
  ...     &  ...      & ... &  ...      & ...\\
z_{n+1,1} & z_{n+1,2} & ... & z_{n+1,n} & 1  \\
\end{array}
\right] \right)^2=
\frac{(-1)^{n+1}}{2^{n}} \cdot \Delta(c_1,...,c_{n+1}).
$$
\vskip 0.3truecm
\par
\noindent
{\bf Proposition~3} (cf. \cite{Blumenthal}, \cite{Borsuk}).
For any points $c_{1},...,c_{n+k} \in {\C}^n$ ($k=2,3,4,...$) their
Cayley-Menger determinant equals $0$ i.e. $\Delta(c_1,...,c_{n+k})=0$.
\vskip 0.3truecm
\par
\noindent
{\it Proof.} Assume that
$c_{1}=(z_{1,1},...,z_{1,n}),...,c_{n+k}=(z_{n+k,1},...,z_{n+k,n})$.
The points $\wt{c}_{1}=(z_{1,1},...,z_{1,n},0,...,0)$,
$\wt{c}_{2}=(z_{2,1},...,z_{2,n},0,...,0)$,...,
$\wt{c}_{n+k}=(z_{n+k,1},...,z_{n+k,n},0,...,0) \in
{\C}^{n+k-1}$ are affinely dependent.
Since $\varphi_n(c_i,c_j)=\varphi_{n+k-1}(\wt{c}_{i},\wt{c}_{j})$
$(1 \leq i \leq j \leq n+k)$ the Cayley-Menger determinant
of points $c_{1},...,c_{n+k}$ is equal to the Cayley-Menger determinant
of points $\wt{c}_{1},...,\wt{c}_{n+k}$ which equals $0$ according to
Proposition~2.
\vskip 0.3truecm
\par
\noindent
From Proposition~2 we obtain the following Propositions~4a and 4b.
\vskip 0.3truecm
\par
\noindent
{\bf Proposition~4a.} If $d,e>0$ with $2ne^2 \neq (n-1)d^2$ and if
points $c_0,...,c_n \in {\C}^n$
satisfy $\varphi_n(c_i,c_j)=d^2$ for $1 \leq i<j \leq n$
and $\varphi_n(c_0,c_i)=e^2$ for $1 \leq i \leq n$,
then the points $c_0,...,c_n$ are affinely independent.
\vskip 0.3truecm
\par
\noindent
{\bf Proposition~4b.} If $d>0$, $c_{1}$, $c_{2}$, $c_{3}$ $\in$ ${\C}^2$ and
$\varphi_2(c_1,c_2)=3d^2$, $\varphi_2(c_1,c_3)=11d^2$,
$\varphi_2(c_2,c_3)=25d^2$, then $c_{1}$, $c_{2}$, $c_{3}$ are
affinely independent.
\vskip 0.3truecm
\par
\noindent
{\bf Proposition~5} (cf. \cite{Borsuk} p. 127 in the real case).
If points $c_{0},c_{1},...,c_{n} \in {\C}^n$ are affinely
independent, $x,y \in {\C}^n$ and
$\varphi_n(x,c_{0})~=~\varphi_n(y,c_{0})$,$~$
$\varphi_n(x,c_{1})~=~\varphi_n(y,c_{1})$,...,
$\varphi_n(x,c_{n})=\varphi_n(y,c_{n})$,
then $x=y$.
\vskip 0.3truecm
\par
\noindent
{\it Proof.} Computing we obtain that the vector
$\overrightarrow{xy}:=[s_1,...,s_n]$ is perpendicular
to each of the $n$ linearly independent vectors
$\overrightarrow{c_{0}c_{i}}$ ($i=1,...,n$).
Thus the vector $\overrightarrow{xy}$ is perpendicular
to every linear combination of vectors
$\overrightarrow{c_{0}c_{i}}$ ($i=1,...,n$).
In particular, the vector
$\overrightarrow{xy}=[s_1,...,s_n]$ is perpendicular
to the vector $[\bar{s_1},...,\bar{s_n}]$, where $\bar{s_1},...,\bar{s_n}$
denote numbers conjugate to the numbers $s_1,...,s_n$, respectively.
Therefore
$\overrightarrow{xy}=0$ and the proof is completed.
\vskip 0.3truecm
\par
\noindent
{\bf Proposition~6a.} The set
$\left\{\left(\sqrt{11}/5\right)^k \cdot \left(\sqrt{3}~\right)^l: k,l \in \left\{0,1,2,...\right\}\right\}$
is a dense subset of $\left(0,\infty\right)$.
\vskip 0.3truecm
\par
\noindent
{\it Proof.} It suffices to prove that
$\left\{\ln_{3}\left(\left(\frac{11}{25}\right)^k \cdot 3^l\right): k,l \in \left\{0,1,2,...\right\}\right\}$
is a dense subset of~$\R$. Computing we obtain the set
$\left\{l+k \cdot \ln_{3}\left(\frac{11}{25}\right): k,l \in \left\{0,1,2,...\right\}\right\}$
which is a dense subset of $\R$ due to Kronecker's theorem (\cite{Hardy}),
because $\ln_{3}\left(\frac{11}{25}\right)$ is irrational and negative.
(If there were $a,b \in \left\{1,2,3...\right\}$ with
$\ln_{3}\left(\frac{11}{25}\right)=-\frac{a}{b}$, then we would have
$3^a \cdot 11^b=25^b$, which is impossible.)
\vskip 0.3truecm
\par
\noindent
Analogously we obtain:
\vskip 0.3truecm
\par
\noindent
{\bf Proposition~6b.} If $n \in \left\{3,4,5,...\right\}$ then
$\left\{\left(\sqrt{2+2/n}~\right)^k \cdot \left(2/n\right)^l: k,l \in \left\{0,1,2,...\right\}\right\}$
is a dense subset of $\left(0,\infty\right)$.
\vskip 0.3truecm
\par
\noindent
{\bf Proposition~7.} For each $n \in \left\{3,4,5,...\right\}$ there exists
$k(n) \in \left\{0,1,2,...\right\}$ such that
$1/2 \leq \left(2/n\right) \cdot \left(\sqrt{2+2/n}~\right)^{k(n)}<1$.
\vskip 0.3truecm
\par
\noindent
{\it Proof.} It is easy to check that
\begin{displaymath}
k(n):=\min \left\{m \in \left\{0,1,2,...\right\}: 1/2 \leq \left(2/n\right) \cdot \left(\sqrt{2+2/n}~\right)^m \right\}
\end{displaymath}
satisfies our condition.
\vskip 0.2truecm
\par
\noindent
As a basic tool for our considerations we shall use the following lemma:
\vskip 0.3truecm
\par
\noindent
{\bf Lemma~1.} Let $n \in \left\{2,3,4,...\right\}$ and
$d,e \in D_n$ such that $2ne^2 >(n-1)d^2$,
and define $r:=\sqrt{4e^2-2 \frac{n-1}{n}d^2}$.
If there exists $\varepsilon \in D_n$
with $\varepsilon \leq 2r$, then $r \in D_n$.
\vskip 0.3truecm
\par
\noindent
{\it Proof.} We may assume that $\varepsilon \neq r$.
Let $x,y \in {\R}^n$, $|x-y|=r$
According to \cite{Benz}, p.19,
there exist points $p_{1},...,p_{n},\wt{y}$,
$\wt{p}_{1},...,\wt{p}_{n} \in {\R}^n$ such that:
\newpage
\par
\noindent
$|x-p_{i}|=|y-p_{i}|=e$ ($1 \leq i \leq n$),
\\
$|p_{i}-p_{j}|=d$ ($1 \leq i<j \leq n$),
\\
$|x-\wt{y}|=r$,
\\
$|y-\wt{y}|=\varepsilon$,
\\
$|x-\wt{p}_{i}|=|\wt{y}-\wt{p}_{i}|=e$ ($1 \leq i \leq n$),
\\
$|\wt{p}_{i}-\wt{p}_{j}|=d$ ($1 \leq i<j \leq n$).
\vskip 0.3truecm
\par
\noindent
Let
$$
S_{xy}:=S_{y\wt{y}}
\cup
\bigcup_{i=1}^{n} S_{xp_{i}}
\cup
\bigcup_{i=1}^{n} S_{yp_{i}}
\cup
\bigcup_{1\le i<j \leq n} S_{p_{i}p_{j}}
\cup
\bigcup_{i=1}^{n} S_{x\wt{p}_{i}}
\cup
\bigcup_{i=1}^{n} S_{\wt{y}\wt{p}_{i}}
\cup
\bigcup_{1\le i<j \leq n} S_{\wt{p}_{i}\wt{p}_{j}}
$$
and $f: S_{xy} \rightarrow {\C}^n$
preserves unit distance. Since
$$
S_{xy} \supseteq
\bigcup_{i=1}^{n}S_{xp_{i}}
\cup
\bigcup_ {i=1}^{n}S_{yp_{i}}
\cup
\bigcup_{1 \leq i<j \leq n} S_{p_{i}p_{j}}
$$
\par
\noindent
we conclude that $f$ preserves the distances between
$x$ and $p_{i}$ ($1\leq i \leq n$), $y$ and $p_{i}$ ($1\leq i\leq n$),
and all distances between $p_{i}$ and $p_{j}$ ($1\leq i<j\leq n$).
Hence for all $1 \leq i \leq n$ $\varphi_n(f(x),f(p_i))=\varphi_n(f(y),f(p_i))=e^2$
and for all $1 \leq i<j \leq n$ $\varphi_n(f(p_i),f(p_j))=d^2$.
Since $S_{xy} \supseteq S_{y\wt{y}}$ we conclude that
$\varphi_n(f(y),f(\wt{y}))=\varepsilon^2$.
By Proposition~3 the Cayley-Menger determinant
$\Delta(f(x),f(p_1),...,f(p_n), f(y))$ equals $0$.
Therefore
$$
\det \left[
\begin{array}{cccccccc}
 0  &  1  &  1  &  1  & ... &  1  &  1  &  1  \\
 1  &  0  & e^2 & e^2 & ... & e^2 & e^2 &  t  \\
 1  & e^2 &  0  & d^2 & ... & d^2 & d^2 & e^2 \\
 1  & e^2 & d^2 &  0  & ... & d^2 & d^2 & e^2 \\
... & ... & ... & ... & ... & ... & ... & ... \\
 1  & e^2 & d^2 & d^2 & ... &  0  & d^2 & e^2 \\
 1  & e^2 & d^2 & d^2 & ... & d^2 &  0  & e^2 \\
 1  &  t  & e^2 & e^2 & ... & e^2 & e^2 &  0  \\
\end{array}
\right]
=0
$$
\par
\noindent
where $t=\varphi_n(f(x),f(y))$.
Computing this determinant we obtain
$$
(-1)^{n-1} \cdot d^{2n-2} \cdot t \cdot
\left(nt+(2n-2)d^2-4ne^2\right)=0.
$$
\newpage
\par
\noindent
Therefore
$$t=\varphi_n(f(x),f(y))=\varphi_n(f(y),f(x))=r^2$$
or $$t=\varphi_n(f(x),f(y))=\varphi_n(f(y),f(x))=0.$$
Analogously we may prove that
$$\varphi_n(f(x),f(\wt{y}))=\varphi_n(f(\wt{y}),f(x))=r^2$$
or $$\varphi_n(f(x),f(\wt{y}))=\varphi_n(f(\wt{y}),f(x))=0.$$
If $t=0$ then the points $f(x)$ and $f(y)$ satisfy:
\\
\\
\centerline{$\varphi_n(f(x),f(x))=0=\varphi_n(f(y),f(x))$,}
\\
\centerline{$\varphi_n(f(x),f(p_1))=e^2=\varphi_n(f(y),f(p_1))$,}
\\
\centerline{$...$}
\\
\centerline{$\varphi_n(f(x),f(p_n))=e^2=\varphi_n(f(y),f(p_n))$.}
\\
\par
\noindent
By Proposition~4a the points
$f(x),f(p_1),...,f(p_n)$ are affinely independent.
Therefore by Proposition~5 $f(x)=f(y)$ and consequently
$$
\varepsilon^2=\varphi_n(f(y),f(\wt{y}))=\varphi_n(f(x),f(\wt{y}))
\in \left\{r^2,~0 \right\}.
$$
Since $\varepsilon^2 \neq r^2$ and $\varepsilon^2 \neq 0$ we conclude
that the case $t=0$ cannot occur. This completes the proof of Lemma~1.
\vskip 0.3truecm
\par
\noindent
As corollaries we obtain:
\vskip 0.3truecm
\par
\noindent
{\bf Lemma~2.}
For $n \in \left\{2,3,4,...\right\}$, $d \in D_n$ and $m \in \left\{0,1,2,...\right\}$,
$\left(\sqrt{2+2/n}~\right)^m \cdot d \in D_n$.
\vskip 0.3truecm
\par
\noindent
{\it Proof.} For $e=\varepsilon=d$, Lemma 1 yields
$r=\sqrt{2+2/n} \cdot d \in D_n$.
This implies the assertion.
\vskip 0.3truecm
\par
\noindent
{\bf Lemma~3.} Let $a,b \in D_2$ such that
$a \leq \frac{4}{\sqrt{5}}b$.
Then, $\sqrt{4b^2-a^2} \in D_2$.
\vskip 0.3truecm
\par
\noindent
{\it Proof.} Since $a<2b$, we have $r:=\sqrt{4b^2-a^2}>0$.
Now, for $\varepsilon:=a$, from Lemma~1 we get the conclusion.
\newpage
\par
\noindent
{\bf Lemma~4.} For $d \in D_2$ and $l \in \left\{0,1,2,...\right\}$,
$\left(\sqrt{3}\right)^l \cdot d \in D_2$. Moreover,
$\sqrt{11} \cdot d \in D_2$ and $5 \cdot d \in D_2$.
\vskip 0.3truecm
\par
\noindent
{\it Proof.} The first assertion follows from Lemma~2
and the second from Lemma~3 for
$(a,b):=(d,~\sqrt{3} \cdot d)$ and $(a,b):=(\sqrt{11} \cdot d,~3 \cdot d)$, respectively.
\vskip 0.3truecm
\par
\noindent
Additionally, for the plane case we need:
\vskip 0.3truecm
\par
\noindent
{\bf Lemma~5.} If $d \in D_2$ then $\left(\sqrt{11}/5\right) \cdot d \in D_2$.
\vskip 0.2truecm
\par
\noindent
{\it Proof.}
Let $d\in D_2$, $x,y \in {\R}^2$, $|x-y|=\left(\sqrt{11}/5\right) \cdot d$.
Using the notation of Figure~1 we show that
$$
S_{xy}:=S_{y\wt{y}} \cup \bigcup_{i=1}^2 S_{xp_i} \cup \bigcup_{i=1}^2 S_{yp_i} \cup S_{p_1p_2}
\cup
\bigcup_{i=1}^2 S_{x\wt{p_i}}
\cup
\bigcup_{i=1}^2 S_{\wt{y}\wt{p_i}}
\cup S_{\wt{p_1}\wt{p_2}}
$$
where the sets corresponding to distances $\sqrt{3} \cdot d$, $\sqrt{11} \cdot d$
and $5 \cdot d$ are known to exist by Lemma~4.
\vskip 0.3truecm
\par
\centerline{
\beginpicture
\normalsize
\setcoordinatesystem units <1.00cm, 1.00cm>
\setplotarea x from -6 to 6, y from -3 to 5
\setplotsymbol({ .})
\plot
-0.872  0.000
 4.564  3.859
 0.000 -1.000
-3.047 -2.776
-0.872  0.000
-5.436  4.859
 0.000  1.000
 2.175 -1.776 
-0.872  0.000 /
\plot
 0.000  1.000
 0.000 -1.000 /
\plot
 2.175 -1.776
-5.436  4.859 /
\plot
-3.047 -2.776
 4.564  3.859 /
\put{$x$} at -1.072 0
\put{$y$} at 0 -1.40
\put{$\wt{y}$} at 0 1.40
\put{$p_1$} at 4.964 3.859
\put{$p_2$} at -3.247 -2.776
\put{$\wt{p_1}$} at -5.636 4.859
\put{$\wt{p_2}$} at 2.575 -1.776
\endpicture}
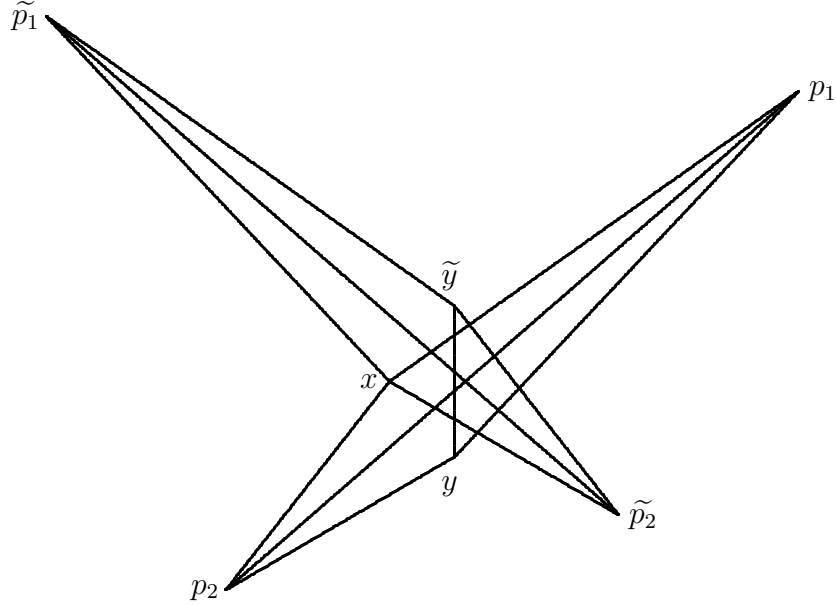
\centerline{Figure 1}
\centerline{$|x-y|=|x-\wt{y}|=\left(\sqrt{11}/5\right) \cdot d$}
\centerline{$|y-\wt{y}|=d$}
\centerline{$|x-p_1|=|y-p_1|=|x-\wt{p_1}|=|\wt{y}-\wt{p_1}|=\sqrt{11} \cdot d$}
\centerline{$|x-p_2|=|y-p_2|=|x-\wt{p_2}|=|\wt{y}-\wt{p_2}|=\sqrt{3} \cdot d$}
\centerline{$|p_1-p_2|=|\wt{p_1}-\wt{p_2}|=5 \cdot d$}
\vskip 0.3truecm
\par
\noindent
Assume that $f: S_{xy} \rightarrow {\C}^2$
preserves unit distance. Since
$$
S_{xy} \supseteq
S_{y\wt{y}}
\cup
\bigcup_{i=1}^{2}S_{xp_{i}}
\cup
\bigcup_ {i=1}^{2}S_{yp_{i}}
\cup
S_{p_1p_2}
$$
\noindent
we conclude that $f$ preserves the distances between
$y$ and $\wt{y}$,
$x$ and $p_{i}$ $(i=1,2)$, $y$ and $p_i$ $(i=1,2)$,
$p_1$ and $p_2$.
By Proposition~3, $\Delta(f(x)$, $f(p_1)$, $f(p_2)$, $f(y))$ equals $0$.
Thus
$$
\det \left[
\begin{array}{ccccc}
0 &  1    &  1    &   1   &   1  \\
1 &  0    & 11d^2 &  3d^2 &   t  \\
1 & 11d^2 &  0    & 25d^2 & 11d^2\\
1 &  3d^2 & 25d^2 &   0   &  3d^2\\
1 &  t    & 11d^2 &  3d^2 &   0  \\
\end{array}
\right]
=0
$$
\vskip 0.3truecm
\par
\noindent
where $t=\varphi_2(f(x),f(y))$.
Computing this determinant we obtain
$$2d^2 \cdot t \cdot \left(11d^2-25t\right)=0.$$
Therefore
$$t=\varphi_2(f(x),f(y))=\varphi_2(f(y),f(x))=\left(\left(\sqrt{11}/5\right) \cdot d \right)^2$$
or $$t=\varphi_2(f(x),f(y))=\varphi_2(f(y),f(x))=0.$$
Analogously we may prove that
$$\varphi_2(f(x),f(\wt{y}))=\varphi_2(f(\wt{y}),f(x))=
\left(\left(\sqrt{11}/5\right) \cdot d\right)^2$$
or $$\varphi_2(f(x),f(\wt{y}))=\varphi_2(f(\wt{y}),f(x))=0.$$
\par
\noindent
If $t=0$ then the points $f(x)$ and $f(y)$ satisfy:
\\
\\
\centerline{$\varphi_2(f(x),f(x))=0=\varphi_2(f(y),f(x))$,}
\\
\centerline{$\varphi_2(f(x),f(p_1))=11d^2=\varphi_2(f(y),f(p_1))$,}
\\
\centerline{$\varphi_2(f(x),f(p_2))=3d^2=\varphi_2(f(y),f(p_2))$.}
\\
\par
\noindent
By Proposition~4b the points
$f(x),f(p_1),f(p_2)$ are affinely independent.
Therefore by Proposition~5 $f(x)=f(y)$ and consequently
$$
d^2=\varphi_2(f(y),f(\wt{y}))=\varphi_2(f(x),f(\wt{y}))
\in \left\{\left(\left(\sqrt{11}/5\right) \cdot d\right)^2,~0 \right\}.
$$
Since $d^2 \neq \left(\left(\sqrt{11}/5\right) \cdot d\right)^2$
and $d^2 \neq 0$ we conclude that
the case $t=0$ cannot occur.
This completes the proof of Lemma~5.
\vskip 0.3truecm
\par
\noindent
Now we obtain:
\vskip 0.3truecm
\par
\noindent
{\bf Theorem~2a.} If $x,y \in {\R}^2$ and
$|x-y|=\left(\sqrt{11}/5\right)^k \cdot \left(\sqrt{3}~\right)^l$
($k,l$ are non-negative integers), then
there exists a finite set $S_{xy}$ with
$\left\{x,y\right\} \subseteq S_{xy} \subseteq {\R}^2$
such that each unit-distance preserving mapping from
$S_{xy}$ to ${\C}^2$ preserves the distance between $x$ and~$y$.
\vskip 0.4 truecm
\par
\noindent
{\it Proof.} Since $1 \in D_2$, Lemmas~4 and~5 imply that
$$
\left\{\left(\sqrt{11}/5\right)^k \cdot \left(\sqrt{3}~\right)^l:
k,l \in \left\{0,1,2,...\right\}\right\} \subseteq D_2.
$$
\vskip 0.3truecm
\par
\noindent
{\bf Theorem~2b.} If $x,y \in {\R}^n$ ($n\geq 3$) and
$|x-y|=\left(\sqrt{2+2/n}~\right)^k \cdot \left(2/n\right)^l$
($k,l$ are non-negative integers),
then there exists a finite set $S_{xy}$ with
$\left\{x,y\right\} \subseteq S_{xy} \subseteq {\R}^n$
such that each unit-distance preserving mapping from $S_{xy}$ to ${\C}^n$
preserves the distance between $x$ and $y$.
\vskip 0.3truecm
\par
\noindent
{\it Proof.}
By Lemma~2 if $d \in D_n$ ($n \geq 2$) then all distances
$\left(\sqrt{2+2/n}~\right)^m \cdot d$ ($m=0,1,2,...$)
belong to $D_n$. By Proposition~7 for each $n \in \left\{3,4,5,...\right\}$
there exists $k(n) \in \left\{0,1,2,...\right\}$ such that
$ 1/2 \leq \rho(n):= \left(2/n\right) \cdot \left(\sqrt{2+2/n}~\right)^{k(n)} < 1$.
For $\varepsilon \in D_n$ we have
$e:=\left(\sqrt{2+2/n}~\right)^{k(n)} \cdot \varepsilon \in D_n$
and $d:=\sqrt{2+2/n} \cdot e \in D_n$.
Therefore, by Lemma~1, if $\varepsilon \in D_n$, then
$\rho(n) \cdot \varepsilon \in D_n$.
Since $1 \in D_{n}$ we conclude that all distances $\rho(n)^m$ ($m=0,1,2,...$)
belong to $D_n$. For each $e \in D_n$ there exists
$m \in \left\{0,1,2,...\right\}$ such that $\varepsilon:=\rho(n)^m \leq \left(4/n\right) \cdot e$.
Applying Lemma~1 for such $\varepsilon$ and for
$d:=\sqrt{2+2/n} \cdot e$ we get:
\vskip 0.3truecm
\par
\noindent
{\boldmath $\left(\ast\right)$}~If $n \in \left\{3,4,5,...\right\}$ and $e \in D_n$,
then $\left(2/n\right) \cdot e \in D_n$.
\vskip 0.3truecm
\par
\noindent
From Lemma~2 and {\boldmath $\left(\ast\right)$} we obtain that for each
non-negative integers $k,l$ we have
$\left(\sqrt{2+2/n}~\right)^k \cdot \left(2/n\right)^l \in D_n$.
This completes the proof of Theorem~2b.
\vskip 0.3truecm
\par
\noindent
{\bf Remark.} If $n \in \{2,3,4,...\}$ and we define the set $D_n$
for $f:S_{xy} \to {\R}^n$, then $D_n$ is equal to the set of all
positive algebraic numbers, see \cite{Aequationes} and \cite{Yerevan}.
\vskip 0.3truecm
\par
\noindent
Finally we state:
\vskip 0.3truecm
\par
\noindent
{\bf Theorem~3.} Each mapping $f:{\R}^n \rightarrow {\C}^n$
($n \geq 2)$
preserving unit distance preserves all distances i.e. satisfies
$$\forall x,y \in {\R}^n ~~\varphi_n(f(x),f(y))=\varphi_n(x,y)=|x-y|^2,$$
provided that $f$ is continuous with respect to the product topologies on
${\R}^n$ and ${\C}^n$.
\vskip 0.3truecm
\par
\noindent
{\it Proof.} Assume that $x,y \in {\R}^n$, $x \neq y$.
Since $D_n \subseteq A_n$:
\begin{description}
\item{{\bf (1)}}
each map $f:{\R}^n \rightarrow {\C}^n$
preserving unit distance preserves all distances belonging to $D_n$
i.e. for all $x,y \in {\R}^n$ if $|x-y| \in D_n$ then
$\varphi_n(f(x),f(y))=\varphi_n(x,y)=|x-y|^2$.
\end{description}
By Theorem~2a
$$\left\{\left(\sqrt{11}/5\right)^k \cdot \left(\sqrt{3}~\right)^l:
k,l \in \left\{0,1,2,...\right\}\right\} \subseteq D_2,$$
so by Proposition~6a $D_2$ is a dense subset of $(0,\infty)$.
By Theorem~2b for all $n \geq 3$
$$\left\{\left(\sqrt{2+2/n}~\right)^k \cdot \left(2/n\right)^l:
k,l \in \left\{0,1,2,...\right\}\right\} \subseteq D_n,$$
so by Proposition~6b for all $n \geq 3$ $D_n$
is a dense subset of $(0,\infty)$.
Therefore, for all $n \geq 2$:
\begin{description}
\item{{\bf (2)}}
there exists a sequence
$\left\{y_k\right\} \subseteq {\R}^n$ tending to $y$
such that for all $k$ $|x-y_k| \in D_n$.
\end{description}
Since $f$ and $\varphi_n$ are continuous,
{\bf (1)} and {\bf (2)} imply that
\begin{displaymath}
\varphi_n(f(x),f(y))=
\lim_{k\rightarrow\infty} \varphi_n(f(x),f(y_k))=
\lim_{k\rightarrow\infty} \varphi_n(x,y_k)=
\varphi_n(x,y)=|x-y|^2.
\end{displaymath}
\vskip 0.3truecm
\par
\noindent
{\bf Acknowledgement.}
The author wishes to thank the anonymous referee for the simplified
proof of Theorem~2b and especially for her or his version
of Theorem~2a with the proof. The author's original version of
Theorem~2a states an analogical result for the distances
$\left(2\sqrt{2}/3\right)^k\cdot \left(\sqrt{3}\right)^l$
with non-negative integers $k,l$
but the proof has been more complicated.
The author also wishes to thank Professor Mowaffaq Hajja
for improvement of the proof of Proposition 1.

A. Tyszka\\
Hugo Ko\l{}\l{}\k{a}taj University\\
Technical Faculty\\
Balicka 104\\
PL-30-149 Krak\'ow\\
Poland\\
e-mail: rttyszka@cyf-kr.edu.pl
\end{document}